\documentclass[final]{siamltex}

\usepackage{amsmath} 
\usepackage{amssymb}

\usepackage[margin=3cm]{geometry} 

\usepackage[all]{xy}

\usepackage[breaklinks,bookmarksopen,bookmarksnumbered]{hyperref}

\usepackage[english]{babel}
\usepackage[utf8]{inputenc}

\usepackage{enumitem}

\usepackage{stackrel} 
\usepackage{tikz-cd} 

\usepackage{ dsfont }

\usepackage{pgfplots}
\pgfplotsset{compat=1.15}
\usepackage{mathrsfs}
\usetikzlibrary{arrows}

\usepackage{algorithmic,algorithm,mathtools}

\newtheorem{teo}{Theorem}[section]

\newtheorem{prop}[teo]{Proposition}

\newtheorem{defin}[teo]{Definition}
\newtheorem{rmk}[teo]{Remark}

\newcommand{\mc}{\mathcal}

\newcommand*{\QED}[1][$\square$]{%
\leavevmode\unskip\penalty9999 \hbox{}\nobreak\hfill
    \quad\hbox{#1}%
}

\def\R{\mathbb{R}}

\def\N{\mathbb N}

\def\P{\mathcal P}
\def\cc{\mathcal C}

\def\cc{\mathcal{C}}

\def\H{\mc{H}}

\newcommand{\m}{\mbox}

\newcommand{\cor}{\textit}

\newcommand{\fine}{\QED\newline}

\DeclareMathOperator{\id}{id}

\DeclareMathOperator{\Proj}{Proj}

\DeclareMathOperator{\var}{var}
\DeclareMathOperator{\KL}{KL}
\DeclareMathOperator{\argmin}{argmin}

\DeclareMathOperator{\sh}{sh}

\def\de{\partial}

\newcommand{\norm}[1]{\left\lVert#1\right\rVert}

\title{A framework for stereo vision via optimal transport}
\author{Mattia Galeotti, Alessandro Sarti, Giovanni Citti}

\begin{document}
\maketitle

\begin{abstract}
We present a theoretical framework for a stereo vision method via optimal transport tools.
We consider two aligned optical systems and we develop the matching between the two
pictures line by line. By considering a regularized version of the optimal transport, we can
speed up the computation and obtain a very accurate disparity function evaluation. Moreover,
via this same method we can approach successfully the case of images with occluded regions.
\end{abstract}

\section{Intro}
Stereo vision is the ability of reconstructing a 3-D representation of a given situation
from the data of two 2-D pictures taken with a small shift. 
This process is clearly fundamental in human vision, has been analyzed through
a great variety of mathematical tools and it is central also in many technological applications. 
In particular stereo vision
techniques employ two cameras with the same focal length looking at the same scene.
The fundamental step in the stereo reconstruction consists in building a correspondence between the two pictures
of the same scene, that is finding for any point on the left image an associated point on the right one.
For any such left-right pair of points, this allows, knowing the relative position of the cameras,
to evaluate the depth in space of the associated real point. One of the main problem with 
the implementation of the reconstruction is the presence
of obstructed regions, that is significant points or objects that are visible from
one optical system but not from the other
because of the obstruction of a nearer object. In \cite{gly95}, \cite{xyw08} and \cite{hka13}
the occlusion problem is widely introduced and reviewed.

Stereo correspondence algorithms can be grouped into those producing sparse
output, usually feature based, and those giving a dense result.
According to the categorization proposed by the important stereo vision taxonomy \cite{schasze02},
dense algorithms can be classified as local and global ones. The local ones,
or window-based, are lean to speed in the accuracy-speed tradeoff. Global ones
are usually slower but more accurate.

 An up-to-date various review of stereo vision techniques
 as well as quick
introduction and brief summary to the state-of-the-art, can be found in \cite{szgmls17}.
In \cite{tlla15} a method for resource-limited systems of stereo vision
algorithms.
Other important review references in stereo vision methods and for the evaluation of stereo
vision algorithms are \cite{bbh03}, \cite{lsg08} and \cite{yiam22}.\newline

In this work, we present the general theoretical framework for a stereo vision process based on an optimal transport
correspondence between the two pictures,
therefore it is a global method with a dense result.
 We consider two cameras looking at the scene
having parallel axis and such that the height of any point is the same from
the two optical systems. This reduces the problem of building the left-right correspondence
to a one-dimensional problem to be carried out on each horizontal section of the pictures,
allowing a gain in efficiency.
Moreover, the reconstruction becomes totally encoded on a single disparity function.


We consider the two (black and white) pictures as functions $I_0,I_1\colon D\to [0,1]$
where $D\subset \R^2$ is a planar domain and we denote the restrictions
to a certain height $y$ by  $I_0^{(y)},I_1^{(y)}$. At first
we treat the case with no obstructed regions, therefore the two horizontal functions
can be seen as finite measures $\nu_0,\nu_1$ over $D^{(y)}=D\cap \R\times \{y\}$ with the same mass, and we
can use the toolbox of optimal transport. 

After normalization, $\nu_0,\nu_1$ are probabilities over $D^{(y)}$. We use the notation
$\pi^{(0)},\pi^{(1)}\colon D^{(y)}\times D^{(y)}\to D^{(y)}$ for the natural
projections, and we call admissible (or admissible plan)
 a probability $\gamma\in\P\left(D^{(y)}\times \P(D^{(y)}\right)$ whose marginals are
$$\pi^{(0)}_\#\gamma=\nu_0,\ \ \pi^{(1)}_\#\gamma=\nu_1.$$
We denote by $\Pi(\nu_0,\nu_1)$ the set of admissible plans
between $\nu_0$ and $\nu_1$.
The problem of optimal transport, by the Kantorovich's formulation,
consists in minimizing 
\begin{equation}
\langle c,\gamma\rangle :=\int c(x,y)d\gamma,
\end{equation}
with $\gamma$ varying in $\Pi(\nu_0,\nu_1)$, where
$c\colon D^{(y)}\times D^{(y)}\to \R$ is called the cost function. In our work
$c$ will be the square difference. It is known (see for example \cite{amgi13})
that under very general assumptions the optimal transport problem has a unique
solution and it induces a map from (a subset of the) left domain to the right domain.

Passing to the discrete setting, anyway, the evaluation of the correspondence above
is very inefficient. For this reason we consider a regularized version of the optimal transport,
we minimize 
\begin{equation}\label{regot}
\langle c,\gamma\rangle -\varepsilon \cdot h(\gamma),
\end{equation}
where $\varepsilon$ is a coefficient and $h$ is the entropy function (see \S\ref{entroreg}).
For any $\varepsilon$ this problem has a unique solution $\gamma^\varepsilon$ in a discrete
setting, and for $\varepsilon$ converging to $0$ we have $\gamma^\varepsilon\to \gamma^\star$,
where $\gamma^\star$ is the maximal entropy probability among the ones minimizing $\langle c,\gamma\rangle$.
The problem of  regularized optimal transport
has been widely treated, and we use \cite{peycut19,cdps17,leo12} as main references.

Minimizing (\ref{regot}) over the space $\cc$ of admissible plans, is the same
of finding the point of $\cc$ minimizing the Kullback–Leibler divergence
from the Gibbs kernel $K^\varepsilon$.
Finding the projection  $\Proj_\Pi^{\KL}(K^\varepsilon)$   
can be done with the so called Sinkhorn's algorithm that we treat in \S\ref{KLsec}, and it has a very good implementation cost.

With respect to our problem, given any probability $\gamma\in \Pi$ we introduce
a function $f_\gamma$ (see Definition~\ref{dispf})
that coincides with the disparity function when $\gamma$ is map induced. In fact $f_\gamma$
works as a good approximation of the disparity function, and our main result
is Theorem \ref{teoconv} where we prove that the Sinkhorn's algorithm
has an exponential rate of convergence with respect to the functions
$f_\gamma$ produced at each step.

We also show that this procedure works very well for artificial non-occluded
cartoons (see Figure~\ref{figspat}), allowing the precise spatial reconstruction with a low implementation cost.\newline

\begin{figure}[h]
\includegraphics[width=14cm]{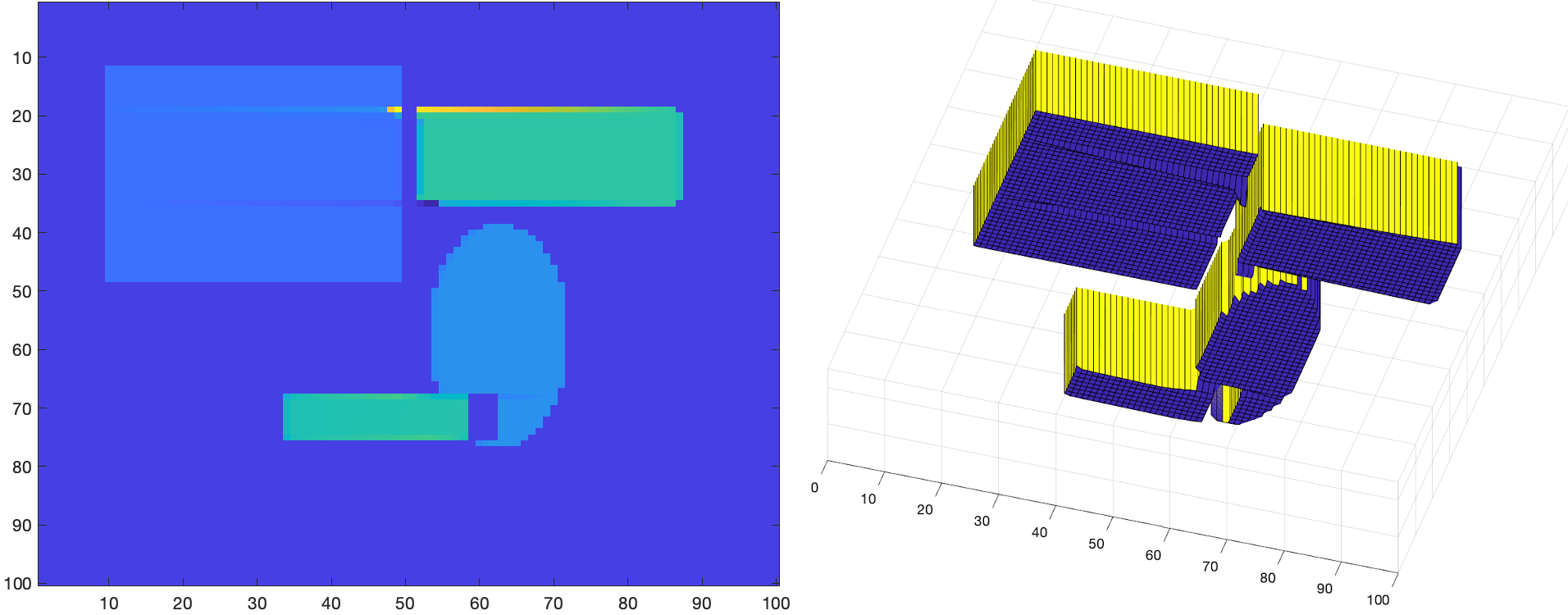}
\caption{We show the disparity function, where different colors correspond
to different disparity values, and the dark blue to points with no mass, therefore including
also occluded regions. On the right we have the three dimensional reconstruction.}\label{figocc}
\end{figure}

In the last section we also try a first approach to the 
case where there are occluded regions. This causes
the masses of the measures $\nu_0,\nu_1$
to be different, therefore the usual theoretical tools of optimal
transport does not work, but we can anyway implement
the Sinkhorn's algorithm.
We show in Theorem \ref{teosep} and Remark \ref{rmksh} that the induced functions $f_\gamma$
converge anyway and the convergence function allows to find the occluded
regions in some simple cases. We also build in Figure \ref{figocc} the reconstruction of 
one of these cases, an artificial cartoon with two object superposing.\newline

\section{Preliminaries on stereo vision}
We recall the basic geometrical intuition 
behind stereo vision and introduce the notations that we will
use along this work.\newline

\subsection{Epipolar geometry}\label{epi}
Every camera, and a human eye, can be modelized as a projection on a 2-dimensional
plane of a 3-dimensional space. More precisely, in the case of the human eye, the projecting surface
is a spherical sector of the retina, but in our discussion we will work with planar projections
because they are more fit to treat camera vision.

Stereo vision is the process of matching points between two planar projections of the same space,
and the consequent reconstruction of the 3-dimensional scene. The general parametrization of the space
through coordinates on two projecting planes is called epipolar geometry (for a wide treatment, see for example \cite{harzis03}).

\begin{figure}[h]
\includegraphics[width=7cm]{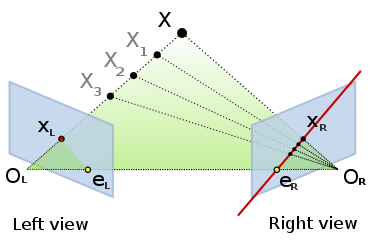}
\caption{The point $X$ in space is projected on the two planes. The red line on the right
is one of the two epipolar lines.}
\end{figure}

We call optical centers the two points $O_L,O_R$ that are the centers for the projections
of any point $X$ on the two planes. We denote by $e_L$ and $e_R$ the projections of the optical centers
$O_R$ and $O_L$ on the left and right plane respectively. For any point $X$ in the space,
we denote by
$X_L,X_R$ its projections on the two planes, and we call \cor{epipolar lines} the lines passing through $X_L,e_L$ on
the left plane, and $X_R,e_R$ on the right plane. The line $O_LO_R$ through the two optical
centers is called \cor{optical line}.\newline

\subsection{Transport problem between two images}
We denote by $I_0$ and $I_1$ the images from the
left and the right optical centers respectively. For both images we have
$$I_j\colon D\to [0,1],\ j=0,1,$$
where the images are represented
by an intensity function on the planar domain $D\subset \R^2$. Therefore
in this approach we are considering images determined by a single
intensity function, such as black and white images.\newline


For any fixed $y$ we consider the restricted domain $D^{(y)}:=D\cap (\R\times \{y\})$
and the restricted images
$$I_j^{(y)}=\left.I_j\right|_{\R\times \{y\}}\colon D^{(y)}\to [0,1],\ j=0,1.$$
We consider $I_0^{(y)}$ and $I_1^{(y)}$ as
the densities of two measures over $D^{(y)}$ and we denote by $\nu_0^{(y)},\ \nu_1^{(y)}$, or simply $\nu_0,\nu_1$, the induced
measure themselves.

We will at first consider the case where $\nu_0$ and $\nu_1$ have the same
total (finite) mass, therefore after normalization they 
are probability measures over the domain $D^{(y)}$. In \S\ref{resu}
we will approach the case where $\nu_0$ and $\nu_1$ have different total mass.
We point out that by the nature of the problem the difference in mass is small with respect to the total mass.\newline

If $\nu_0,\nu_1$ are probability measures
over $D^{(y)}$, our goal is to find
 the disparity shift between the two pictures
by applying the optimal transport tools, 
that is by searching the optimal plan that realizes the infimum (\ref{kantorovich}).
Because of the nature of the problem, we are interested in optimal
plans that are induced by transport maps. In fact, the solution to the optimal transport problem
is map-induced in this case.

\begin{prop}\label{propJ}
Consider a closed interval $J\subset \R$ and two probability measures $\nu_0,\nu_1\in \P(J)$
that are absolutely continuous with respect to the classic Lebesgue measure over $J$.
If we consider the quadratic cost $c\colon J\times J\to \R$ such that $c(x,y)=(x-y)^2$,
then there exists a unique optimal plan $\gamma$ in $\Pi(\nu_0,\nu_1)$
and it is induced by a map $T\colon J\to J$. 
\end{prop}
\proof If the measures are absolutely continuous with respect to the Lebesgue
measure, then they are regular and we conclude by \cite[Theorem 1.26]{amgi13}.\fine

\begin{rmk}\label{rmkJ}
We observe that this is precisely our case, where $J=D^{(y)}$ and
$\nu_0,\nu_1$ are induced by the densities $I_0^{(y)},I_1^{(y)}$.
Indeed, $I_0^{(y)},I_1^{(y)}$ are continuous almost everywhere and therefore
the two induced measures are absolutely continuous with respect
to the classic Lebesgue measure.\newline
\end{rmk}

\subsection{A model of the human optical system}\label{hos}
In our work, we will suppose that the two projection planes are both parallel to the optical line,
and this line is considered to be horizontal, in such a way that
the vertical coordinate of a same point, is the same in the two projections.

As both image planes and the optical line 
 are parallel, we can consider the distance $b$
of the image plane (also called projection plane) from the optical system as a single invariant of our system.

For any point $X$ in space,  we denote by $x_L$ and $x_R$ the (horizontal) coordinates
of its two projections. 
We are interested in finding the disparity between the position
of the object seen from an optical center, and the position from the other one.
We suppose that objects ``at infinity'' are seen without disparity, that is
we fix the coordinates on the projecting planes in such a way that parallel lines
from the optical centers intercept the projecting plan at the same coordinate $x_L=x_R$.


\begin{figure}[h]
\definecolor{wrwrwr}{rgb}{0.3803921568627451,0.3803921568627451,0.3803921568627451}
\definecolor{rvwvcq}{rgb}{0.08235294117647059,0.396078431372549,0.7529411764705882}
\begin{tikzpicture}[line cap=round,line join=round,>=triangle 45,x=1.0cm,y=1.0cm]
\clip(-7.,-1.5) rectangle (5.,4.);
\draw [line width=2.pt,color=wrwrwr] (-2.97434,3.52613)-- (-1.5571985824775605,-0.86);
\draw [line width=2.pt,color=wrwrwr] (-2.97434,3.52613)-- (1.7816756770739022,-0.86);
\draw (-5.5,1.5) node[anchor=north west] {projection plane};
\draw (-5,-0.3) node[anchor=north west] {optical line};
\draw [line width=2.pt,dash pattern=on 1pt off 1pt,color=wrwrwr] (1.2388742595514624,0.82)-- (1.7816756770739022,-0.86);
\draw [line width=2.pt,color=wrwrwr] (-2.84396243131991,0.82)-- (3.3660237810885145,0.82);
\draw [line width=2.pt,color=wrwrwr] (-3.083319501605459,-0.86)-- (3.9777140718182515,-0.86);
\begin{scriptsize}
\draw [fill=rvwvcq] (-2.1,0.82) circle (2.5pt);
\draw[color=rvwvcq] (-1.8,1.1) node {$X_L$};
\draw [fill=rvwvcq] (-0.04,0.82) circle (2.5pt);
\draw[color=rvwvcq] (0.10145929469393566,1.1) node {$X_R$};
\draw [fill=rvwvcq] (-2.97434,3.52613) circle (2.5pt);
\draw[color=rvwvcq] (-2.883855276367502,3.769053961544046) node {$X$};
\draw [fill=wrwrwr] (1.7816756770739022,-0.86) circle (2.0pt);
\draw[color=wrwrwr] (2,-1.2) node {$O_R$};
\draw [fill=wrwrwr] (-1.5571985824775605,-0.86) circle (2.0pt);
\draw[color=wrwrwr] (-1.8,-1.2) node {$O_L$};
\draw [fill=wrwrwr] (1.2388742595514624,0.82) circle (2.0pt);
\draw[color=wrwrwr] (1.4578160263120497,1.1) node {$X_L'$};
\end{scriptsize}
\end{tikzpicture}
\caption{The object $X$ is projected to $X_L$ and $X_R$ respectively from the optical center $O_L$ and $O_R$.
The coordinate $x_R$ on the right is set such that $x_r(X_L')=0$,
where $O_RX_L'$ is a segment parallel to $XO_L$. }\label{figgeo}
\end{figure}
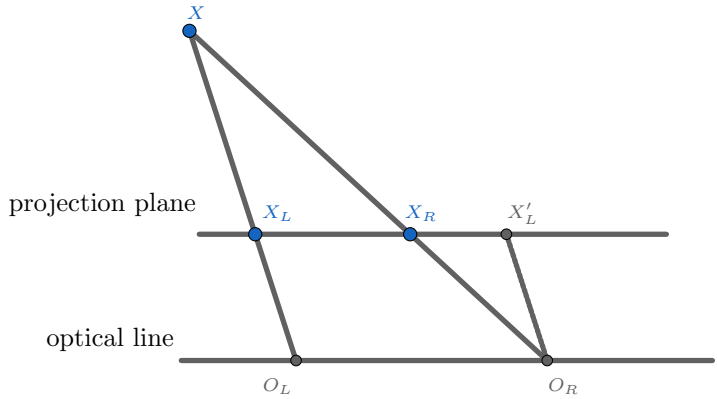

As it appear from the image, for an object at distance $b$ from the 
optical line, if $\ell$ is the distance of $X$ from the optical line
and $\ell_0$ is the distance $|O_LO_R|$ between the optical centers, then
\begin{equation}\label{shift}
x_L-x_R=|X_RX_L'|=\frac{b}{\ell}\cdot \ell_0.
\end{equation}

We remark that there is only one coordinate considered because we are supposing that
the optical system is aligned with the vertical ax. This is compatible with the ability
of human vision to distinguish `up' from `down'. The matching procedure by optimal transport between the 
left and right projections
will be applied at every
$y$ coordinate separately. This speeds up the efficiency of our algorithm, because
it scales with $\lambda^2$ if $\lambda$ is the length of the horizontal window that we consider. Instead, 
 if we had considered optimal transport on a $2$-dimensional domain, the scaling factor would have been $\lambda^4$.\newline

\section{Optimal transport and stereo vision}\label{optstereo}

\subsection{Definition of the optimal transport problem}
We recall briefly the formulations of the optimal transport problem.
Our main reference is Ambrosio-Gigli guide \cite{amgi13}.
We recall that a Polish space is a complete and separable metric space.
Consider two Polish spaces $X,Y$, and a Borel \cor{cost} function $c\colon X\times Y\to \R_{\geq0}\cup \{+\infty\}$.
In the Monge definition of  the optimal transport problem, we start from two probability
measures $\nu_0$ and $\nu_1$ in $X$ and $Y$ respectively.  
We are interested in finding (if it exists) a function $T\colon X\to Y$
sending $\nu_0$ in $\nu_1$ and such that 
 the total transport cost is minimized. We recall that the measure $T_{\#}\nu_0$ over $Y$
 is defined by
 $$T_\#\nu_0(B)=\nu_0(T^{-1}(B))\ \forall B\in \mc B(Y).$$
 Therefore, the optimal transport problem consists in achieving
 \begin{equation}\label{monge}
 \inf_{T_{\#}\nu_0=\nu_1}\int_Xc(x,T(x))d\mu,
 \end{equation}
where the infimum is taken over all the functions $T$ respecting the constraint $T_{\#}\nu_0=\nu_1$.\newline

Kantorovich generalizes this approach, by considering transport plans.
\begin{defin}
An \cor{admissible transport plan} between $\nu_0$ and $\nu_1$ is a measure $\gamma$
over $X\times Y$
such that $\pi^X_\#\gamma=\nu_0$ and $\pi^Y_\#\gamma=\nu_1$, or equivalently
\begin{align*}
\gamma(A\times Y)&=\nu_0(A)\ \ \forall A\in\mc{B}(X)\\
\gamma(X\times B)&=\nu_1(B)\ \ \forall B\in\mc{B}(Y).
\end{align*}
We denote the set of admissible transport plans between $\nu_0$ and $\nu_1$ by
$\Pi(\nu_0,\nu_1)$.\newline
\end{defin}

The optimal transport problem becomes the detection of
\begin{equation}\label{kantorovich}
\inf_{\gamma\in \Pi(\nu_0,\nu_1)}\int_{X\times Y}c(x,y)d\gamma.
\end{equation}
\begin{defin}
An optimal plan is an admissible plan achieving the infimum of (\ref{kantorovich}).
\end{defin}
\begin{rmk}
We introduce a more compact notation for the coupling on the right. If $c$
is a cost function and $\gamma\in \P(X\times Y)$, we define
$$\langle c,\gamma\rangle :=\int_{X\times Y}c(x,y)d\gamma.$$
\end{rmk}

In fact,  problem (\ref{kantorovich}) has a solution under very general conditions.

\begin{teo}[see {\cite[Theorem 4.1]{villa08}}]\label{teoex}
Consider two measures $\nu_0$, $\nu_1$ in $\P(X)$ and $\P(Y)$ respectively.  
If the cost function $c$ is lower semicontinuous and bounded from below, then there exists
an optimal plan $\gamma$ for the functional
$\gamma\mapsto \langle c,\gamma\rangle$,
among all $\gamma\in \Pi(\nu_0,\nu_1)$.\newline
\end{teo}

We are interested in the case where an optimal plan is induced by a transport map $T$,
which means that the optimal plan $\gamma$ respects $(\id\times T)_{\#}\nu_0=\gamma$, or equivalently
the solution of (\ref{monge}) exists and coincides with the more general solution of (\ref{kantorovich}).
We saw in Proposition \ref{propJ} and Remark \ref{rmkJ}
that in our cases the opportune conditions are verified for the optimal
plan to be induced by a transport map. The same
is true in the discrete setting that we approach from now on.

\begin{rmk}\label{discrete}
When we work in a discrete setting, we consider the domain of
the probability measures $\nu_0,\nu_1$ to be $X=Y=[1,d]\subset \N$.
Therefore, $\nu_0,\nu_1$ are horizontal vectors in $[0,1]^d$.
For the sake of a clearer notation we denote by $i$ the coordinate on $X$
and by $j$ the coordinate on $Y$.

An admissible transport plan $\gamma\in \Pi(\nu_0,\nu_1)$ is a $d\times d$ matrix with
positive entries such that 
\begin{equation}\label{constr}
\mathds{1}\gamma=\mu,\ \ \mathds{1}\gamma^T=\nu,
\end{equation}
where $\mathds{1}=(1,1,\dots,1)$ is the horizontal vector whose coordinates are all $1$.
In this presentation, the value $\gamma_{ij}$ is the mass transported from $i$ to $j$
for any $i,j=1,\dots,d$. Furthermore, given a cost function 
$c\colon[1,d]^2\to \R_{\geq0}$, the cost coupling becomes
$$\langle c,\gamma\rangle :=\sum_{i,j=1}^d c(i,j)\cdot \gamma_{ij}.$$\newline
\end{rmk}

\subsection{Entropic regularization}\label{entroreg}
We rename the solution of the
optimal transport problem
$$L(\mu,\nu):=\inf_{\gamma\in \Pi(\mu,\nu)}\langle c,\gamma\rangle.$$
By Theorem \ref{teoex} this has 
 a solution if we work with a quadratic cost $c$. Moreover, the solution is unique in the case of measures induced
by densities over a closed interval of $\R$, as in the case we are considering
in treating stereo vision.

For a given admissible plan $\gamma$, we denote by $h(\gamma)$
its entropy
$$h(\gamma):=-\int_{X\times Y}(\log(\gamma)-1)d\gamma.$$
We consider the optimal transport problem after an entropic regularization,
and denote by $L^\varepsilon$ the regularized infimum,
\begin{equation}\label{regular}
L^\varepsilon(\mu,\nu):=\inf_{\gamma\in\Pi(\mu,\nu)}\left(\langle c,\gamma\rangle-\varepsilon \cdot h(\gamma)\right).
\end{equation}

For a wide introduction to the entropic regularized problem see \cite{cdps17}, \cite[\S4]{peycut19} and \cite{leo12}, here
we recall the main aspects. 

\begin{prop}   
Given $\nu_0,\nu_1\in \P(H)$ for a closed interval $H\subset \R$, and absolutely continuous with respect
to the Lebesgue measure, 
there exists unique an admissible plan $\gamma^\varepsilon$ achieving the infimum~(\ref{regular}),
and also
$$L^\varepsilon (\nu_0,\nu_1)\xrightarrow{\varepsilon\to 0}L(\nu_0,\nu_1).$$
\end{prop}
\begin{rmk}
For a proof of the above result see \cite[Theorem 2.7]{cdps17}, where is also proved
that if $\gamma$ is the unique optimal plan for $L(\nu_0,\nu_1)$, then
$\gamma^\varepsilon$ narrowly converges to $\gamma$.\newline
\end{rmk}

In the discrete setting $h(\gamma)=-\sum_{i,j}\gamma_{ij}(\log(\gamma_{ij})-1)$
and the result above is also true (see \cite[Proposition 4.3]{peycut19}).
In particular, if for a given cost function $c\colon [1,d]^2\to \R_{\geq0}$
we consider the Gibbs kernel
\begin{equation}\label{gibbs}
K^{\varepsilon}\in \R^{d\times d}\ \m{s.t.}\ K_{ij}^{\varepsilon}:=e^{-\frac{c(i,j)}{\varepsilon}}\ \ \ \forall i,j=1,\dots,d,
\end{equation}
then the unique solution of (\ref{regular}) has the form
\begin{equation}\label{solu}
\gamma_{ij}^\varepsilon=u_i^\varepsilon K_{ij}^\varepsilon v_j^\varepsilon,
\end{equation}
for some vectors $u^\varepsilon,v^\varepsilon\in \R^d$. The proof 
is straightforward after introducing the Lagrangian
of (\ref{regular}),
$$\mc{L}(\gamma,f,g):=\langle c,\gamma\rangle -\varepsilon h(\gamma)-\langle f,\mathds{1}\gamma-\mu\rangle-\langle g,\mathds{1}\gamma^T-\nu\rangle.$$

\begin{rmk}\label{hreg}
In the discrete setting the optimal plan $\gamma$ is not necessarily unique.
Anyway, by strict concavity of $h$ there exists unique $\gamma^\star\in \Pi(\nu_0,\nu_1)$
realizing the minimum $L(\nu_0,\nu_1)$ and maximizing $h(\gamma^\star)$.
We have (see \cite[Proposition 4.1]{peycut19}) that $\gamma^\varepsilon \xrightarrow{\varepsilon\to 0}\gamma^\star$.\newline
\end{rmk}

\section{Reconstructing disparity shifts with optimal transport}
In this section we find an equivalent formulation of the regularized transport
problem in terms of a known divergence function. Then, we build an associated disparity function
to any transport plan.

\subsection{Projecting via the Kullback–Leibler divergence}\label{KLsec}
Divergences are functions that share some properties with distances
but are not symmetric nor they respect the triangular inequality.
They are frequently used to compare measures and probability distributions.

\begin{defin}
If $\P$ is a differentiable manifold, a divergence on $\P$ is a function $F\colon \P\times \P\to \R$
verifying the following  properties:
\begin{enumerate}
\item $F(x,y)\geq 0$ for all $x,y\in \P$;
\item $F(x,y)=0$ if and only if $x=y$;
\end{enumerate}
\end{defin}

One of the most used divergences in probability and statistics is the Kullback–Leibler divergence.

\begin{defin}
The Kullback–Leibler divergence, or KL divergence,
is defined for any probability measures $\gamma$, $\alpha$ over a 
domain $Z$. If $\gamma$
is absolutely continuous with respect to $\alpha$, then
$$\KL(\gamma|\alpha):=\int_Z\log\left(\frac{d\gamma}{d\alpha}\right)d\gamma,$$
otherwise $\KL(\gamma|\alpha)=+\infty$.
\end{defin}

\begin{rmk}
In a discrete setting as the one introduced above, $\gamma$ is always absolutely continuous
with respect to $\alpha$ and the KL divergence becomes
$$\KL(\gamma|\alpha)=\sum_{i,j}\gamma_{ij}\cdot \log\left(\frac{\gamma_{ij}}{\alpha_{ij}}\right).$$
\end{rmk}

The KL divergence behaves similarly to the regularized cost if we consider the Gibbs kernel $K^\varepsilon$
introduced in (\ref{gibbs}). Indeed,
$$\varepsilon\cdot \left(\KL(\gamma|K^\varepsilon)-1\right)=\langle c,\gamma\rangle -\varepsilon \cdot h(\gamma),$$
therefore minimizing the right side of the above equation for $\gamma\in \Pi(\nu_0,\nu_1)$,
is equivalent to minimize the KL divergence with $K^\varepsilon$ for $\gamma$ in the same space.

\begin{rmk}
Observe that $K^\varepsilon$ is not a probability distribution. It is however possible to define the KL divergence,
although it does not have all the divergence properties above listed. 
Therefore, the optimal transport regularized problem (\ref{regular}) becomes the problem
of finding the projection with respect to $\KL$ of the Gibbs kernel $K^\varepsilon$ on $\Pi(\mu,\nu)$.
\end{rmk}

We introduce two affine subspaces of the probability space $\P([1,d])$
where the conditions (\ref{constr}) are verified.

\begin{defin}
For any $\mu,\rho\in \P([1,d])$ we denote by $\cc_\mu^1$ the subspace of $\P([1,d]^2)$
of probability measures $\gamma$ verifying $\mathds{1}\gamma=\mu$.
Analogously, we denote by $\cc_\nu^2$ the subspace of $\gamma$ verifying $\mathds{1}\gamma^T=\rho$.\newline
\end{defin}

\begin{rmk}\label{rmkproj}
Observe that $\Pi(\nu_0,\nu_1)=\cc_{\nu_0}^1\cap \cc_{\nu_1}^2$. 
In order to state an algorithmic implementation
for the resolution of the regularized optimal transport, we consider alternatively the projections
on $\cc_{\nu_0}^1$ and~$\cc_{\nu_1}^2$,
\begin{align*}
\gamma^{(2k+1)}&:=\Proj^{\KL}_{\cc_{\nu_0}^1}\left(\gamma^{(2k)}\right),\\
\gamma^{(2k+2)}&:=\Proj^{\KL}_{\cc_{\nu_1}^2}\left(\gamma^{(2k+1)}\right),
\end{align*}
for any $k\in \N$, after initializing the procedure with $\gamma^{(0)}:=K^\varepsilon$.
\end{rmk}

\begin{prop}
For any function $K$ over $[1,d]^2$ and any probability $\mu\in \P([1,d])$,
the plan $\overline \gamma\in \cc_\mu^1$ minimizing $\KL(\overline \gamma|K)$
is
$$\overline\gamma_{ij}:=\overline u_iK_{ij},$$
where the vector $\overline u\in \R^d$ is the unique vector
such that $\overline \gamma$ respects $\mathds{1}\overline\gamma=\mu$.

Analogously the same is true for $\cc_\rho^2\ni\overline\gamma:=(K_{ij}\overline v_j)$
where $\overline v$ is the unique vector such that $\overline \gamma$
respects $\mathds{1}\overline\gamma^T=\rho$.
\end{prop}
\proof  We prove the first result. The Lagrangian of $\inf_{\gamma\in \cc_\mu^1}\KL(\gamma|K)$ is
$$\mc L(\overline\gamma,f)=\sum_{i,j}(\overline\gamma_{ij}\log(\overline\gamma_{ij})-\overline\gamma_{ij}\log(K_{ij}))-\langle \beta,\mathds{1}\overline\gamma-\mu\rangle,$$
where $\beta\in \R^d$ is a variable $d$-dimensional vector.
If we derive in $\overline\gamma_{ij}$ we get
$$\frac{\de \mc L}{\de \overline\gamma_{ij}}=\log(\overline\gamma_{ij})+1-\log(K_{ij})-\beta_i.$$
If we impose $\frac{\de \mc L}{\de \overline\gamma_{ij}}=0$ for any $i,j\in [1,d]$,
then we must have
$$\overline\gamma_{ij}=K_{ij}\cdot e^{\beta_i-1}.$$
We define $\overline u_i:=e^{\beta_i-1}$ and we get the matrix form in the theorem.
The condition $\mathds{1}\overline\gamma=\mu$ is necessary, but by construction
is now also sufficient and gives the unicity of the solution.\fine

This motivates the well known Sinkhorn's algorithm.
After initializing the
procedure with the vectors $u^{(0)}=v^{(0)}=\mathds{1}_d$, we have
\begin{align*}
u^{(k+1)}&:=\nu_0\oslash (K^\varepsilon \cdot v^{(k)}),\\
v^{(k+1)}&:=\nu_1\oslash ((K^\varepsilon)^T\cdot u^{(k+1)}),
\end{align*}
where the division on the right sides is the entry-wise division.
And with this notation we have $\left(\gamma_{ij}^{(2k+1)}\right)=\left(u^{(k+1)}_iK^\varepsilon_{ij}v^{(k)}_j\right)$
and $\left(\gamma_{ij}^{(2k)}\right)=\left(u^{(k)}_iK^\varepsilon_{ij}v^{(k)}_j\right)$.
As $\cc_{\nu_0}^1$ and $\cc_{\nu_1}^2$ are affines, by \cite{breg67}
we know that the iterative process converges to $\Proj_{\Pi(\nu_0,\nu_1)}^{\KL}\left(K^\varepsilon\right)$
which is precisely the matrix $\gamma^\varepsilon$ as expressed in~(\ref{solu}).

\begin{rmk}
Observe that the vectors $u^\varepsilon,v^\varepsilon$ depend on the 
initializing vectors $u^{(0)},v^{(0)}$ but the unique optimal plan $\gamma^\varepsilon$
does not.\newline
\end{rmk}


\subsection{Disparity functions}\label{dispfun}
Any stereo vision reconstruction of 3-dimensional scene
is equivalent to the datum of a disparity shift at any point of the right projection $I_1$,
 that is a vector encoding the
shift of that point when passing from the right optical system to the left one.
In our setting, the two optical systems are vertically aligned, therefore
the disparity shift is a scalar positive value giving the horizontal shift (to the right).

Therefore, given a spatial configuration $C$ of some objects
and the associated stereo vision picture pair $I_0,I_1$, at any $y$ coordinate the spatial reconstruction
must be a disparity function $f\colon S\to \R_{\geq0}$, where $J$ is a closed interval in $\R$
and $S\subset J$. 

In this section
we work in the discrete setting, therefore
the domain of the disparity functions is a subset of $J=[1,d]$.
Furthermore, we consider
disparity functions resulting from images where there is no occlusion, meaning
that there are no areas of the depicted objects that are invisible from one of the optical systems
because of the superposition of other objects. We will consider
cases with occlusions in \S\ref{resu}.\newline


We associate to any plan $\gamma\in \P([1,d]^2)$
a function $f_\gamma$ and we will prove below (see Proposition~\ref{f=f})
that it is a non-negative disparity function in the case of stereo vision induced plans.

\begin{defin}\label{dispf}
If $\gamma\in \P([1,d]^2)$, for any $i\in [1,d]$ such that
$\mathds{1}\gamma(i)=\sum_j\gamma_{ij}>0$,
we define the following function,
$$f_\gamma(i):=\frac{\sum_j\gamma_{ij}\cdot j}{\mathds{1}\gamma(i)}-i.$$\newline
\end{defin}

\begin{defin}
We call object any compact closed submanifold $X$ in $\R^3$ whose boundary $\partial X$
is smooth $\mc L_{\partial X}$-a.e.~and $\mc L_{\partial X}$ is the canonical Lebesgue measure
on $\partial X$.
\end{defin}

Because of the aligned disposition of the optical systems,
we consider any horizontal section of the configuration, and its
left and right projection to the optical line.\newline

Consider a stereo vision picture pair $I_0,I_1\colon [1,d]\to [0,1]$ obtained
in the discrete setting from a configuration $C$. If $I_0,I_1$ have the same finite mass at any $y$,
 we can normalize, obtaining $\nu_0^{(y)},\nu_1^{(y)}\in \P([1,d])$. In order
to apply the optimal transport theorems, we consider the cost function
$c\colon[1,d]^2\to \R_{\geq0}$ defined by
$$c(i,j):=(i-j)^2.$$

\begin{lemma}\label{lemver1}
Consider $\gamma$ an optimal plan for the problem $L(\nu_0,\nu_1)=\inf \langle c,\gamma\rangle$,
where the infimum is taken in $\Pi(\nu_0,\nu_1)$, and
$\gamma_{ij}>0$ and $\gamma_{i'j'}>0$ for some $i,i',j,j'\in [1,d]$.
If $i'>i$, then $j'\geq j$.
Equivalently, if $\gamma_{ij}>0$ then there are no $i'>i$ and $j'<j$ such that $\gamma_{i'j'}>0$.
\end{lemma}
\proof  Let's suppose that there exists $i'>i$ and $j'<j$ such that $\gamma_{ij}>0$
and $\gamma_{i'j'}>0$. We want to show that the plan $\gamma$ is not optimal, contradicting the hypothesis.
In order to do this, we introduce 
$$\gamma_{\min}:=\min(\gamma_{ij},\gamma_{i'j'}),$$
and define another plan $\gamma'\in\Pi(\nu_0,\nu_1)$, equal to $\gamma$ except that
\begin{align*}
\gamma_{ij}'&:=\gamma_{ij}-\gamma_{\min}\\
\gamma_{ij'}'&:=\gamma_{ij'}+\gamma_{\min}\\
\gamma_{i'j'}'&:=\gamma_{i'j'}-\gamma_{\min}\\
\gamma_{i'j}'&:=\gamma_{i'j}+\gamma_{\min}
\end{align*}
In order to show that $\langle c,\gamma'\rangle < \langle c,\gamma\rangle$
it suffices to observe that
$$(i'-j')^2+(i-j)^2>(i'-j)^2+(i-j')^2.$$\fine

\begin{defin}
We say that a point $P$ of an object $X$ of the configuration
projects to the left (or the right) optical system if the segment $O_LP$ (or $O_RP$)
does not intersect any object except at $P$.

In this case, we say that the coordinate $x_L(P)$ (or $x_R(P)$)
represents the projected point $P$.
\end{defin}
\begin{defin}\label{occluded}
We call 
occluded regions those subsets of the space that are `visible'
from one of the two optical systems but not from the other.
More precisely
\begin{align*}
H_0&:=\{x\in [1,d]|\ x=x_R(P)\ P\m{ is not represented on the left o.s.}\};\\
H_1&:=\{x\in [1,d]|\ x=x_L(P)\ P\m{ is not represented on the right o.s.}\}.
\end{align*}
\end{defin}

In the example below, we see that areas $H_j$
correspond to elements that are covered by a closer (to the optical system) object,
and therefore are invisible for one eye or the other. For an introduction to the occlusion problem
see for example \cite{gly95, hka13}.


\begin{figure}[h]
\includegraphics[width=12cm]{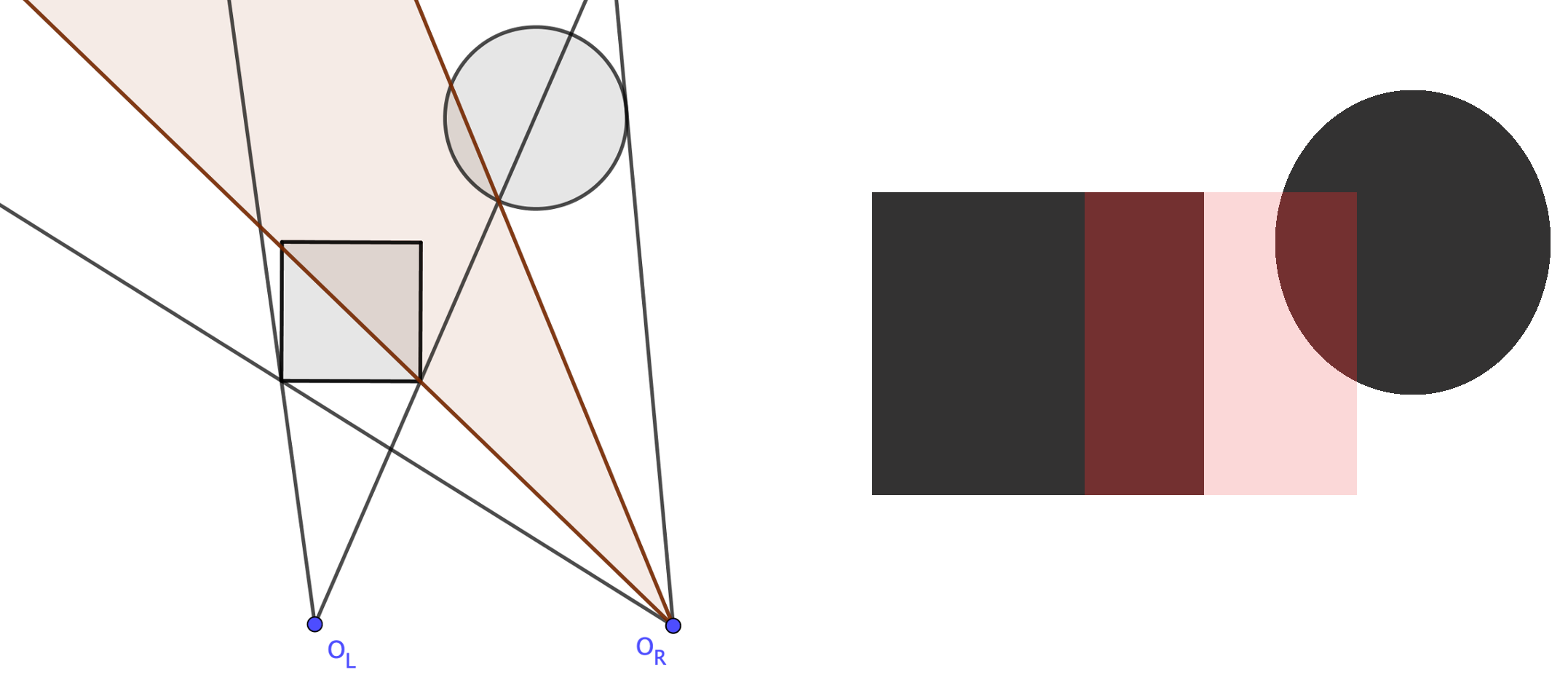}
\caption{On the left there is the disposition in space of the two objects, seen from above. On
the right we have the image from the right optical center, the shadowed area
is $H_0$.}
\end{figure}

We want to define the disparity function associated with a spatial
configuration of a finite number of objects.
We start by considering very simple configurations.
\begin{defin}\label{nocart}
We call non-occluded cartoon a configuration where any object 
is included in a plane or a line parallel to the optical
line, and every point on any object is represented on both optical systems.
\end{defin}

We will drop the `non-occluded' hypothesis in the following section,
adapting our methods to cases with non-empty occluded regions.

Given a non-occluded cartoon $C$, we consider
the horizontal section at a chosen height $y$
and define an associated disparity function $f^{(C,y)}$
that we denote simply $f^{(C)}$ when there is no risk of confusion.
We consider the subset $S\subset [1,d]$
of $x_R$ representing a point $P$ 
on any object of $C$.

\begin{defin}
If $x_L$ and $x_R$ are defined as in \S\ref{hos},
for any $x_R\in S$ representing a point $P$, we have
$$f^{(C)}(x_R):=x_L(P)-x_R(P).$$
\end{defin}

By definition, for any non-occluded cartoon $C$ and any
$y$ coordinate,
$f^{(C,y)}$ is a non-negative function.\newline






Consider a stereo vision picture pair $I_0,I_1\colon [1,d]\to [0,1]$ obtained
in the discrete setting from
a non-occluded cartoon $C$. By construction, $I_0,I_1$ have the same finite mass at any $y$,
and we denote again by $\nu_0,\nu_1\in \P([1,d])$ the probabilities
obtained after the normalization.


\begin{prop}\label{f=f}
Given a non-occluded cartoon $C$ and any optimal plan $\gamma$
between the induced probabilities $\nu_0,\nu_1$ for the cost function $c$ defined above,
we have $f_\gamma=f^{(C)}$ and therefore $f_\gamma$
is a non-negative function.
\end{prop}
\proof  We prove inductively that for any point $i\in S=\{x\in [1,d]|\ \nu_0(x)>0\}$,
if $i$ represents a point $P$ of the configuration, that is $i=x_R(P)$,
then $\gamma_{ij}=\nu_0(i)$ if $j=x_L(P)$, and $\gamma_{ij}=0$ for any other $j$.

This is true for the minimum $i\in [1,d]$ such that $\nu_0(i)>0$,
because if $i=x_R(P)$ then $j=x_L(P)$ is the minimum $j$ such that $\nu_1(j)>0$
and $\nu_1(j)=\nu_0(i)$ by construction. And we conclude with Lemma~\ref{lemver1}.

Suppose now that the proposition is true up to $i_0$, and $i_0=x_R(P)$, $j_0=x_L(P)$.
Consider $i_1:=\argmin \{i>i_0,\ \nu_0(i)>0\}$, and $j_1:=\argmin\{j>j_0,\ \nu_1(j)>0\}$.
By the absence of obstruction we must have $i_1=x_R(Q)$ for some point
$Q$ of the configuration such that $j_1=x_L(Q)$. Moreover, by Lemma \ref{lemver1}
the mass in $j_1$ must be `filled' with the mass `coming from' $i_1$,
therefore $\gamma_{i_1j_1}=\nu_0(i_1)=\nu_1(j_1)$.\fine

\section{Convergence results}
We show the convergence of the estimates for the optimal
plan and also of the associated disparity functions.
Furthermore, we introduce here a slight variation of the problem we considered above by
 that the target measure $\nu_1$ is a probability measure while the total
mass of the source measure is $m(\nu_0)>1$. This will in particular account for
the difference among the two pictures due to the occlusion of some objects
by other objects that are nearer to the optical system. The case where the target measure $\nu_1$
 carries the greater mass is theoretically equivalent to the first case, but with slight differences
 in the applied consequences in stereo vision settings. 

\subsection{Convergence of associated functions}
In order to treat the convergence of Sinkhorn's algorithm, we introduce
for any vector $w\in \R^d$ the variation norm
$$\norm{w}_{\var}:=\max_iw_i-\min_iw_i,$$
and for any two vectors $u,u'\in \R^d_{>0}$ the Hilbert metric 
$$d_\H(u,u'):=\norm{\log (u)-\log(u')}_{\var}.$$
Observe that in order for $d_\H$ to be a metric, we have
to quotient out the relation $u\sim u'$
if $u=\lambda\cdot u'$ for some $\lambda\in \R_{>0}$.
Over $\R^d_{>0}\slash \sim$, $d_\H$ respects the 
triangle inequality and it is a distance.\newline

For every positive matrix $K\in \R_{>0}^{d\times d}$,
we define
$$\eta(K):=\max_{i,j,i',j'}\frac{K_{ij}\cdot K_{i'j'}}{K_{i'j}\cdot K_{ij'}},$$
and moreover
$$\lambda(K):=\frac{\sqrt{\eta(K)}-1}{\sqrt{\eta(K)}+1}.$$
Therefore $0<\lambda(K)<1$ and by the following theorem
$K$ acts as a contraction on the cone of positive vectors.

\begin{teo}
For $K$ as above and $u,u'\in \R^d_{>0}$, 
$$d_\H(Ku, Ku')\leq\lambda(K)\cdot d_\H(u,u').$$
\end{teo}
For a proof see \cite[Theorem 4.1]{peycut19}.
Moreover, as a direct consequence of the theorem
above we can state a convergence rate result for the sequences $u^{(k)}$
and $v^{(k)}$ converging to $u^\varepsilon,v^\varepsilon$
that we treated in~\S\ref{KLsec}.

We observe that if $u^\varepsilon(i)=0$ for some $i\in [1,d]$,
then $\nu_0(i)=0$ and $u^{(k)}(i)=0$ for any $k\geq1$. And the same is true if $v^\varepsilon(j)=0$
for $\nu_1(j)=0$ and $v^{(k)}(j)=0$ if $k\geq1$. Therefore the convergence on null coordinates is completed at the first iteration,
and for the evaluation of the distance between two vectors we can focus only on non-null coordinates.
Thus, up to considering a coordinate subspace $\R^{d'}\subset \R^d$
we can suppose that $d_\H(u^{(k)},u^\varepsilon)$ and $d_\H(v^{(k)},v^\varepsilon)$
are well defined for any $k\in \N_{>0}$, or equivalently that their coordinates
are everywhere non-null.

\begin{lemma}\label{lemD}
If $u^{(k)},v^{(k)}$ are the convergent series of vectors
induced by the Sinkhorn's algorithm as explained above,
we have the following convergence rate
\begin{align*}
d_\H(u^{(k)},u^\varepsilon)&=O(\lambda(K^\varepsilon)^{2k})\\
d_\H(v^{(k)},v^\varepsilon)&=O(\lambda(K^\varepsilon)^{2k}).
\end{align*}
\end{lemma}
\proof For any triple of everywhere non-null vectors $u,u',a\in \R^d_{>0}$, we have
that $d_\H(u,u')=d_\H(a\oslash u,a\oslash u')$.
Moreover, because of the iteration rules of the Sinkhorn's algorithm, we know
that 
\begin{align*}
d_\H(u^{(k)},u^\varepsilon)&=d_\H(\nu_0\oslash(K^\varepsilon v^{(k-1)}),\nu_0\oslash(K^\varepsilon v^\varepsilon))\\
&=d_\H(K^\varepsilon v^{(k-1)},K^\varepsilon v^\varepsilon)\\
&\leq \lambda(K)\cdot d_\H(v^{(k-1)},v^\varepsilon).
\end{align*}
In the same way we can prove $d_\H(v^{(k-1)},v^\varepsilon)\leq \lambda(K^\varepsilon)\cdot d_\H(u^{(k-1)},u^\varepsilon)$,
and this concludes the proof.\fine

We are able to translate the convergence rate result above in a result about the 
convergence of disparity functions. We already considered the unique admissible plan $\gamma^\varepsilon$
 minimizing $L_\varepsilon(\nu_0,\nu_1)$ and the series of plans $\gamma^{(2k+1)}=u^{(k+1)}K^\varepsilon v^{(k)}$ for any $k\in \N_{>0}$,
 obtained via the Sinkhorn's algorithm.
We define 
\begin{align*}
f_\varepsilon&:=f_{\gamma^\varepsilon}\\
f_{(k)}&:=f_{\gamma^{(2k+1)}}\newline
\end{align*}
We show that $f_{(k)}$ converges to $f_\varepsilon$
and estimate the convergence rate.
\begin{teo}\label{teoconv}
For $\varepsilon\to 0$, we have $f_{(k)}\to f_\varepsilon$
and moreover
$$\norm{f_{(k)}-f_\varepsilon}_\infty= O(\lambda(K^\varepsilon)^{2k}).$$
\end{teo}
\proof Observe that for any $i\in[1,d]$,
$$u^{(k+1)}(i)=\frac{\nu_0(i)}{\sum_j K^\varepsilon_{ij}v_j^{(k)}},$$
for any $k\in \N$. Therefore we can rewrite
the disparity function for $\gamma^{(2k+1)}$ evaluated at 
any $i$ gives
$$f_{(k)}(i)+i=\frac{\sum_{j}u_i^{(k+1)}K_{ij}^\varepsilon v_j^{(k)}\cdot j}{\mathds{1}\gamma^{(2k+1)}(i)}=\frac{\sum_{j}K_{ij}^\varepsilon v_j^{(k)}\cdot j}{\sum_j K^\varepsilon_{ij}v^{(k)}_j},$$
because $\nu_0(i)=\mathds{1}\gamma^{(2k+1)}(i)$.
By Lemma \ref{lemD}, if we call $d_\H(v^{(k)},v^\varepsilon)=:A$,
then
$$e^{-A}\cdot\frac{\sum_jK_{ij}^\varepsilon v^\varepsilon_j\cdot j}{\sum_jK_{ij}^\varepsilon v_j^\varepsilon}\leq
\frac{\sum_j K_{ij}^\varepsilon v^{(k)}_j\cdot j}{\sum_j K_{ij}^\varepsilon v^{(k)}_j}\leq
e^A\cdot \frac{\sum_jK_{ij}^\varepsilon v^\varepsilon_j\cdot j}{\sum_jK_{ij}^\varepsilon v_j^\varepsilon}.$$
By straightforward calculations this implies
$$\left |f_{(k)}(i)-f_\varepsilon(i)\right|\leq d\cdot \left(1-e^{-A}\right),$$
and again by Lemma \ref{lemD} this implies the theorem.\fine







\subsection{Shifted projections}\label{shipro}

In the case where $\nu_1$ is a probability measure but $\nu_0$ is not, and it has
a total mass $m_0:=m(\nu_0)>1$, the main problem in evaluating the regularized optimal
plan is that $\Pi(\nu_0,\nu_1)=\cc_{\nu_0}^1\cap\cc_{\nu_1}^2=\varnothing$.
In the following we are going to use the Sinkhorn's procedure of successive projections, 
and show that we will obtain two separate convergent series.\newline

Observe that the measure $\widetilde\nu_0:=\frac{1}{m_0}\cdot\nu_0$
is a probability measure by definition. Moreover, observe that 
any element of $\cc_{\nu_0}^1$ has mass $m_0$ and any element
of $\cc_{\widetilde\nu_0}^1$ is a probability measure.

\begin{defin}
We call shifting $\sh\colon \cc_{\nu_0}^1\to \cc^1_{\widetilde\nu_0}$
the map that sends $\mu\mapsto \frac{1}{m_0}\cdot \mu$. When there is no risk
of confusion we use the notation $\widetilde\mu=\sh(\mu)$. 
\end{defin}
\begin{rmk}
Observe that
this map preserves the affine structure with respect
to the vector space of $0$ mass functions. Moreover, $\sh$
is invertible as a map of affine spaces
\end{rmk}

We work again in the discrete setting where $\nu_0,\nu_1$ are
defined over $[1,d]$.
For any measure $\alpha$ over $[1,d]^2$ consider the two sequences 
\begin{align*}
\gamma_\alpha^{(2k+1)}&:=\Proj_{\cc_{\nu_0}^1}^{\KL} \left(\gamma_\alpha^{(2k)}\right)\\
\gamma_\alpha^{(2k+2)}&:=\Proj_{\cc_{\nu_1}^2}^{\KL} \left(\gamma_\alpha^{(2k+1)}\right)
\end{align*}
for any $k\in \N$, after the initialization $\gamma_\alpha^{(0)}:=\alpha$, analogous to what we did in Remark \ref{rmkproj}
for the Gibbs kernel $K^\varepsilon$.
In this case we obtain two separate convergence results.

\begin{teo}\label{teosep}
In the case above, we have that the even sequence converges
$$\gamma_\alpha^{(2k)}\xrightarrow{k\to+\infty}\left(\Proj_{\Pi(\widetilde\nu_0,\nu_1)}^{\KL} (\alpha)\right)\in \cc_{\nu_1}^2,$$
and the odd sequence converges too, but to another limit
$$\gamma_\alpha^{(2k+1)}\xrightarrow{k\to+\infty}\sh^{-1}\left(\Proj_{\Pi(\widetilde\nu_0,\nu_1)}^{\KL}(\alpha)\right)\in \cc_{\nu_0}^1.$$
\end{teo}
\proof We are going to show that for any measure $\mu \in \cc_{\nu_0}^1$, 
minimizing the KL divergence $\KL(\mu|\alpha)$ is equivalent to minimizing
$\KL(\widetilde\mu|\alpha)$. Indeed,

\begin{align*}
\KL(\mu|\alpha)&= \int \log\left(\frac{ d\mu}{d\alpha}\right)d\mu\\
&=\int\log\left(\frac{m_0\cdot d\widetilde\mu}{d\alpha}\right)m_0\cdot d\widetilde\mu\\
&=\int m_0\cdot \left(\log(m_0)+\log\left(\frac{d\widetilde\mu}{d\alpha}\right)\right)d\widetilde \mu\\
&=m_0\cdot\left(\log(m_0)+\KL(\widetilde\mu|\alpha)\right).
\end{align*}
Moreover, $\tilde \mu$ lies in $\cc^1_{\widetilde \nu_0}$ which intersects
 with $\cc^2_{\nu_1}$ at $\Pi(\widetilde \nu_0,\nu_1)$.
Therefore the sequence of probabilities $\widetilde\gamma_\alpha^{(1)},\gamma_\alpha^{(2)},\widetilde\gamma_\alpha^{(3)},\gamma_\alpha^{(4)},\dots$
respects the conditions of the sequence in Remark \ref{rmkproj} with initialization equal to $\alpha$.
Therefore both $(\widetilde\gamma_\alpha^{(2k+1)})_k$ and $(\gamma_\alpha^{(2k)})_k$
converge to $\Proj^{\KL}_{\Pi(\widetilde\nu_0,\nu_1)}(\alpha)$, proving the theorem.\fine

\begin{rmk}
Observe that the limit of $\gamma_\alpha^{(2k+1)}$ lies
in the subset of $\cc_{\nu_0}^1$ of minimal KL divergence from $\cc_{\nu_1}^2$, 
that is $\sh^{-1}(\Pi(\widetilde\nu_0,\nu_1))$.
\end{rmk}
\begin{rmk}\label{rmksh}
If $\gamma\in \Pi(\nu_0,\nu_1)$ and $\widetilde\gamma=\sh\gamma$, then
the associated functions $f_\gamma$ and $f_{\widetilde\gamma}$ coincide.
Therefore in the succession above, also the disparity function converges,
and the converge rate result of Theorem \ref{teoconv} is still true.\newline
\end{rmk}

\section{Numerical results}\label{resu}
First we consider the case of a pair of probability measures obtained
by taking two artificially drawn images that are shifted one respect to the other and have no obstructions.
Then we also apply the same method for cartoons with a simple obstruction, that is
one object covers another one on the same scene.

\subsection{An example of non-occluded cartoon}
In this section we consider the discrete domain $D=[1,d]\times[1,h]\subset \N^2$. Therefore the
two images are functions 
$$I_j\colon[1,d]\times [1,h]\to [0,1],\ \ j=0,1,$$
and for every $y\in[1,h]$ coordinate, $I_j^{(y)}\colon [1,d]\to [0,1]$
are the intensity functions that determine (after normalization)
the probabilities $\nu_0^{(y)},\nu_1^{(y)}$. We consider the 
case of a non-occluded cartoon as by Definition \ref{nocart}.
This represents a first modelization of a real configuration of objects
where no-object is superposed to any other from both of the optical systems.\newline

By Proposition \ref{f=f}, in the case of a quadratic cost and
a non-occluded cartoon, the only optimal plan $\gamma$ is the one
induced by the natural matching between the two stereo vision pictures.

In order to evaluate this matching plan for every horizontal line,
our procedure consists in implementing
the Sinkhorn's algorithm described in \S\ref{KLsec} to every 
height $y$, because the regularized transport solution is a good
approximation of $\gamma$ if the regularization coefficient
is sufficiently small (see Remark~\ref{hreg}). Other than the $I_j$ functions,
we have to introduce the entropic coefficient $\varepsilon$ and the
number of iterations \texttt{niter}.

\begin{algorithm}[H]
        \begin{algorithmic}
	\caption{Computation of the transport plan between $I_0$ and $I_1$ via the Sinkhorn's method } \label{alg1}
          \STATE \textbf{Input}  $I_0$, $I_1$, $\varepsilon$, \texttt{niter}.
          \STATE Set $d=$ \texttt{size}$(I_0,2)=$ \texttt{size}$(I_1,2)$, and $h=$ \texttt{size}$(I_0,1)=$ \texttt{size}$(I_0,1)$;
          \STATE Define the cost matrix $c\in \R^{d\times d}$ such that $c(i,j):=(i-j)^2$;
          \STATE Define the kernel $K^{\varepsilon}(i,j):=e^{-\frac{c(i,j)}{\varepsilon}}$;
          \STATE Initialize $v^{(0)}:=$ \texttt{ones}$(1,d)$;
          \FOR{$y=1,\dots,h$}
          \STATE Set $\nu_0:=I_0(y,\colon)$ and $\nu_1:=I_1(y,\colon)$;
          \FOR{ $k=1,\dots, \m{\texttt{niter}}$}
           \STATE $u^{(k)}:=\nu_0\oslash(K^\varepsilon v^{(k-1)})$;
           \STATE $v^{(k)}:=\nu_1\oslash((K^\varepsilon)^Tu^{(k)})$;
           \ENDFOR
           \STATE Set $u^\varepsilon:=u^{(\m{\texttt{niter}})}$ and $v^\varepsilon:=v^{(\m{\texttt{niter}})}$;
           \STATE Set $\gamma^\varepsilon(\colon,\colon,y)=\diag(u^\varepsilon)K^\varepsilon\diag(v^\varepsilon)$.
           \ENDFOR
        \end{algorithmic}
      \end{algorithm}    	

We develop this with the cartoon in the image below, taking
three horizontal sections as examples.\newline

\begin{figure}[h]
\includegraphics[width=11cm]{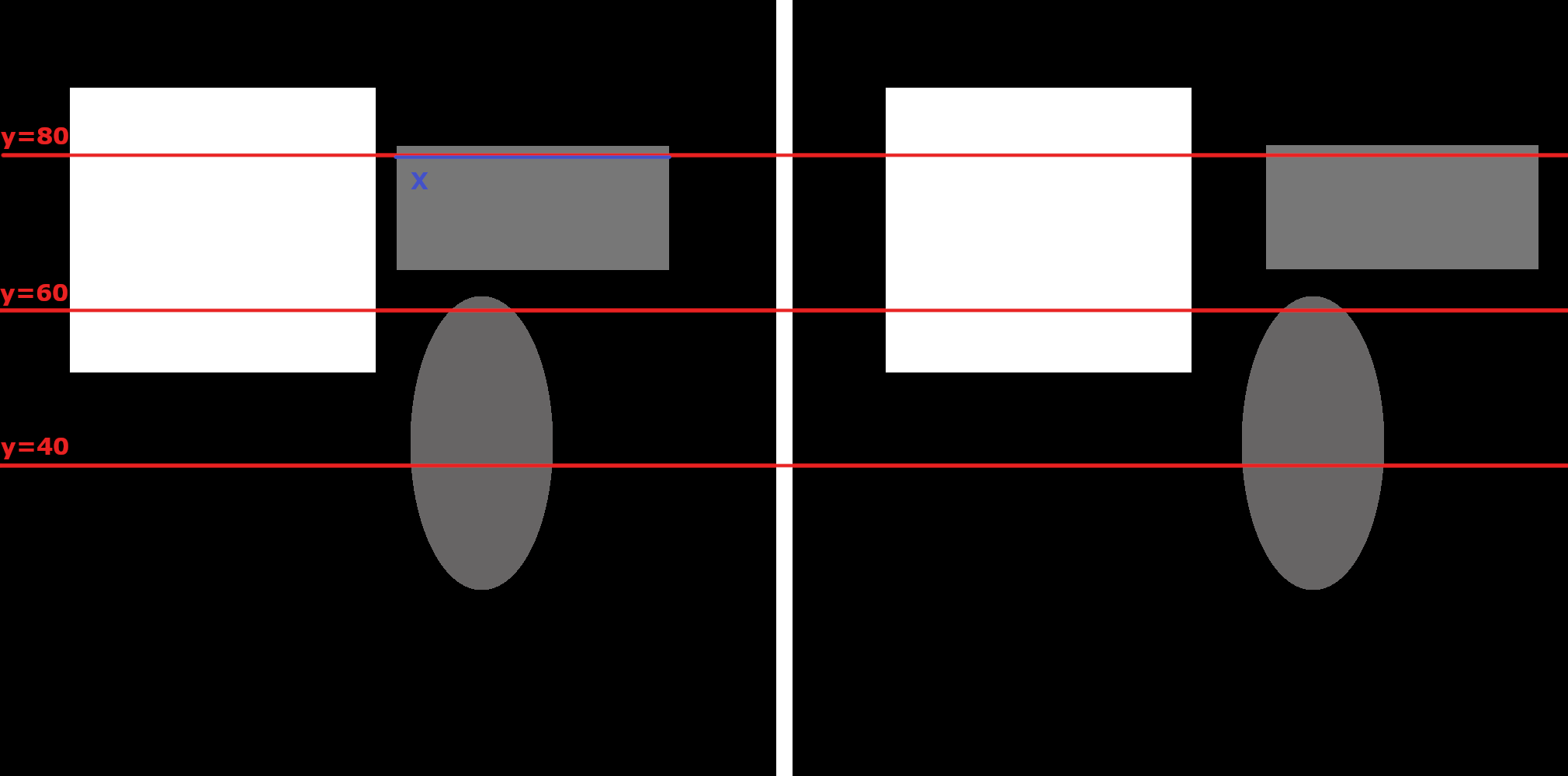}
\caption{Two pictures of the same non-obstructed cartoon taken from two different optical systems. We
marked three horizontals sections for $y=40,60,80$, the rectangular object $X$ and we also
marked in blue the intersection $D^{(80)}\cap X$.}
\end{figure}

\begin{figure}[h]
\includegraphics[width=11cm]{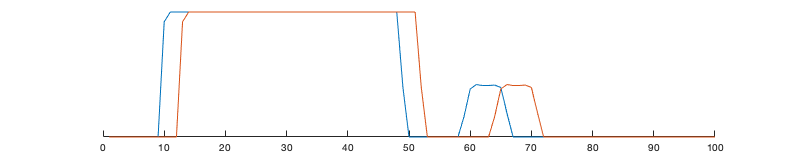}
\caption{The intensity functions $I_0^{(60)}$ (in blue) and $I_1^{(60)}$ (in orange),
at the horizontal section $y=60$.}
\end{figure}

\begin{figure}[h]
\includegraphics[width=15cm]{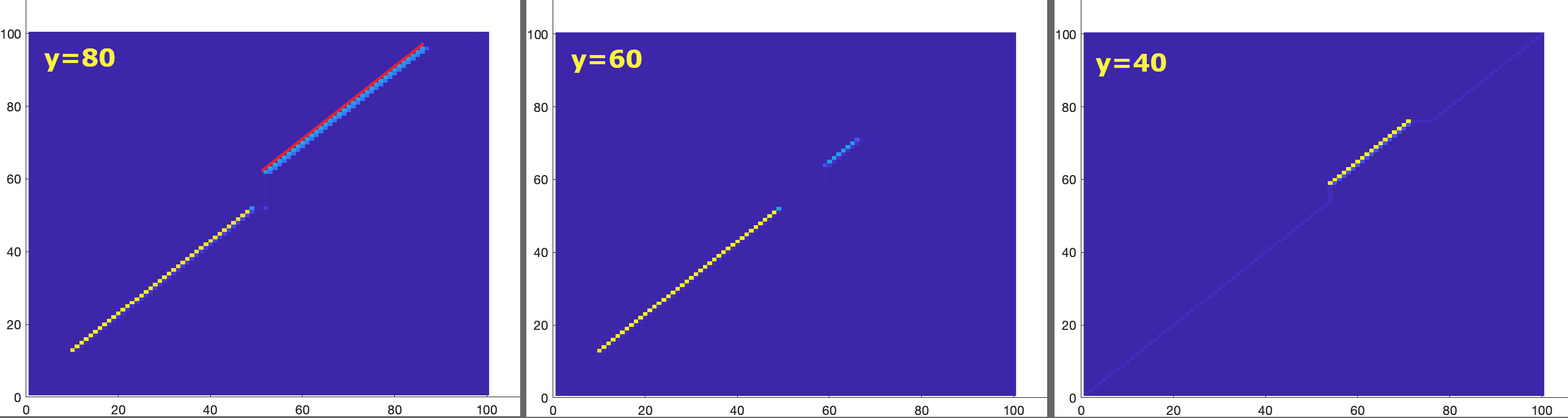}
\caption{An approximation of the optimal transport plan matrices for the marked lines $y=80,60,40$, after a sigmoid passage that emphasizes
the non-null entries. The marked segment on the first image is the intersection $X\cap D^{(80)}$. The matrices
are obtained with the log version of the Sinkhorn's algorithm described below.}\label{imgg}
\end{figure}

Implementing our algorithm with \texttt{niter}$=10^5$ and $\varepsilon=0.1$
we obtain the optimal transport plans and the disparity
functions depicted in Figure \ref{imgg}. The disparity results agree
with the cartoon configuration we used.
Finally we develop a 3-dimensional reconstruction
of the cartoon, that only depends on the disparity functions
and the geometrical invariants of the optical system.\newline

\begin{figure}[h]
\includegraphics[width=15cm]{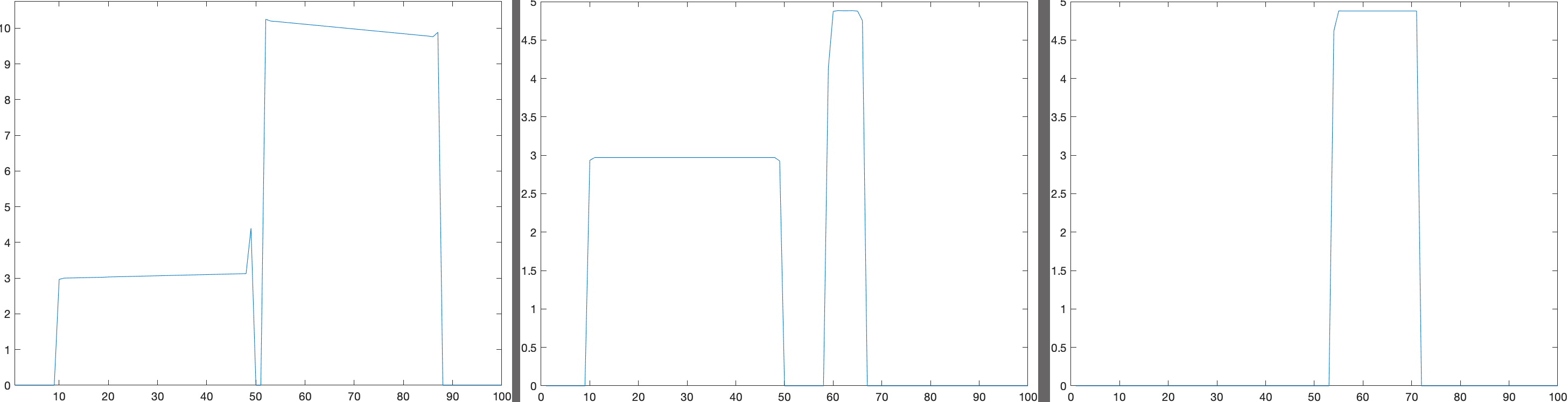}
\caption{The disparity functions obtained at $y=80,60,40$. See Definition \ref{dispf}.}
\end{figure}

\begin{figure}[h]
\includegraphics[width=10cm]{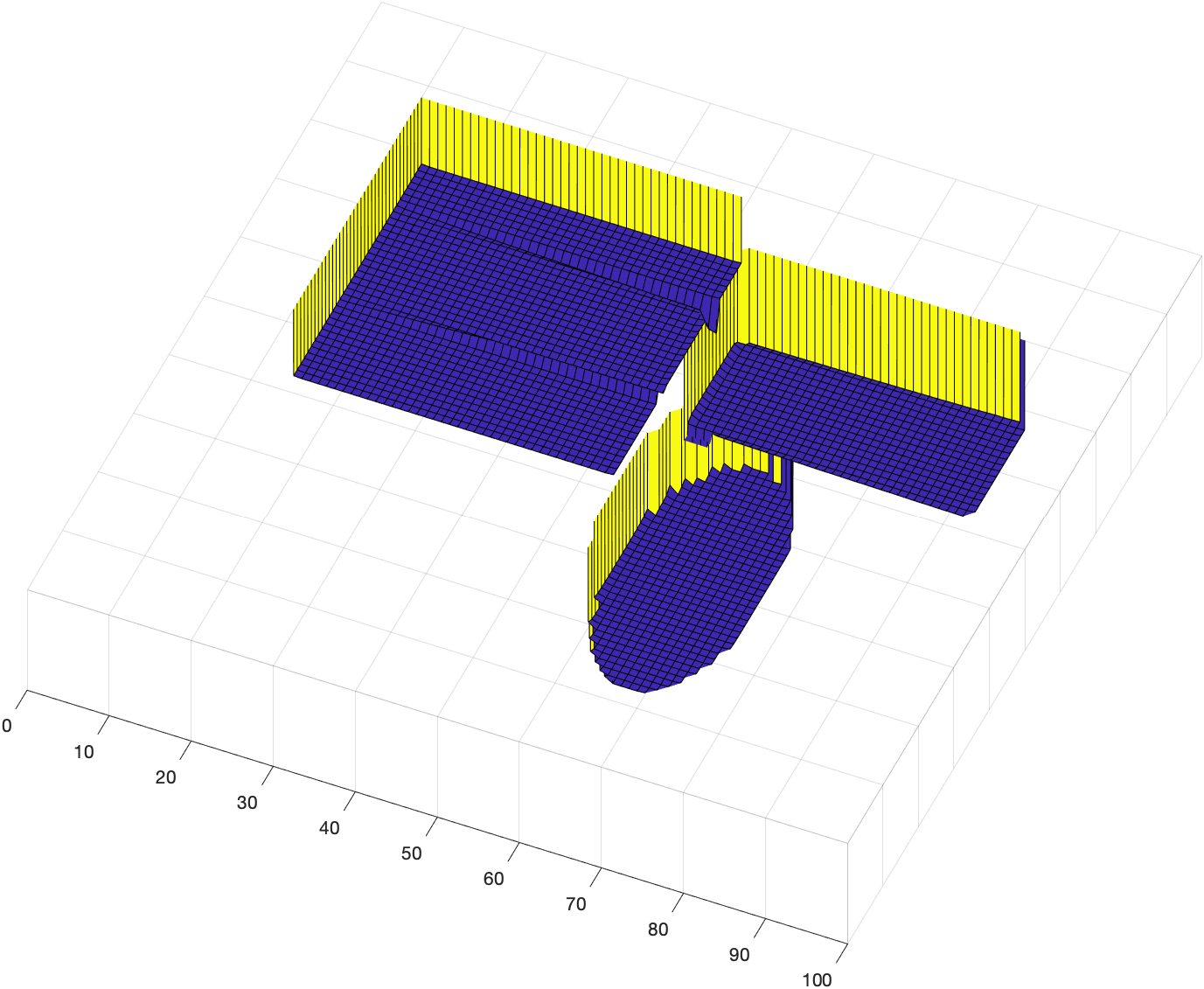}
\caption{A three dimensional reconstruction after the estimate (line by line) of the disparity shift.
We supposed that $\ell_0=10$ and $b=1000$ (see Figure \ref{figgeo}), and
we have a conversion coefficient $\beta=2$ between pixels and the length unit, therefore $\ell=\frac{b\cdot \ell_0}{\beta\cdot (x_L-x_R)}$.}\label{figspat}
\end{figure}

\subsection{From shifted projections to the recovery of occluded regions}
Starting from the results of \S\ref{shipro}, it is possible
to determine the occluded regions of a given configuration
if they exist. We consider in particular the case of a cartoon configuration
where $H_0=\varnothing$ while $H_1\neq \varnothing$, meaning that there is some point
represented in the left optical system that is not represented
in the right optical system, but the inverse is not true. If $H_1\neq \varnothing$
then $m(I_0^{(y)})\geq m(I_1^{(y)})$ for some $y$.

Differently from the non-obstructed case, we obtain  that the two measures $\nu_0^{(y)}$
and $\nu_1^{(y)}$ may have different masses. Therefore, after dividing by $m(I_1^{(y)})$
$\nu_1$ is a probability measure while $m(\nu_0)\geq 1$.
We focus on those $y$ where the last inequality is strict,
and therefore the mass quotient $\varphi^{(y)}:=m(I_1^{(y)})\slash m(I_0^{(y)})$
is strictly lower than $1$. When there is no risk of confusion we denote
this quotient by $\varphi$.\newline

\begin{defin}
We call compression coefficient of a point $x\in[1,d]$
de derivative of $f_\gamma$ if it is well defined in $x$.
\end{defin}

We observe that by the results of \S\ref{dispfun}
if $\nu_0,\nu_1$ are two probabilities then the compression coefficient
is $0$ everywhere it is defined. If $m(\nu_0)>1$ then, 
as stated in Remark \ref{rmksh}, the disparity function
$f_\gamma$ can be obtained by considering the optimal transport
problem between $\widetilde \nu_0$ and $\nu_1$. This means that the compression
coefficient can take non-null values.\newline

We consider the leftmost object in our configuration and denote it by $X$.
As we are matching $\widetilde \nu_0$ instead of $\nu_0$, then the mass of $X$ is
reduced by a $\varphi$ factor. 
Using Lemma \ref{lemver1}, we observe that $\widetilde \nu_0|_X$ `fills'
the corresponding area of $\nu_1$ slower than in the non-obstructed case.
We state without proof that the new compression coefficient
at points in $X$ approximates $\varphi-1$.\newline

If $i_0$ is the coordinate of the leftmost point of $X$, then
$f_\gamma(i_0)$ recovers correctly the shift of $X$. This means that we
can recover also the region occluded by $X$. Indeed, if we know the rightmost point
$i_1$ of $X$, then by checking how much of the mass of $I_1$ lies
at the right of $i_1+f_\gamma(i_0)$ we can recover
the area that is covered by $X$. This means
recovering the point $i_2$ such that the segment
$[i_1+1,i_2]$ is occluded by $X$.
We consider then $\nu_0':=\nu_0-\nu_0|_{[i_1+1,i_2]}$.
If $m(\nu_0')=1$ then
we have recovered all the occluded area; if not,
we repeat the same procedure by eliminating the object $X$ both from $\nu_0$ and $\nu_1$.\newline

In the example below, we considered a case with $H_0=\varnothing$ and $H^1\neq \varnothing$,
in particular focusing at $y=30$ we obtain an occluded region and recover
it.

\begin{figure}[h]
\includegraphics[width=16cm]{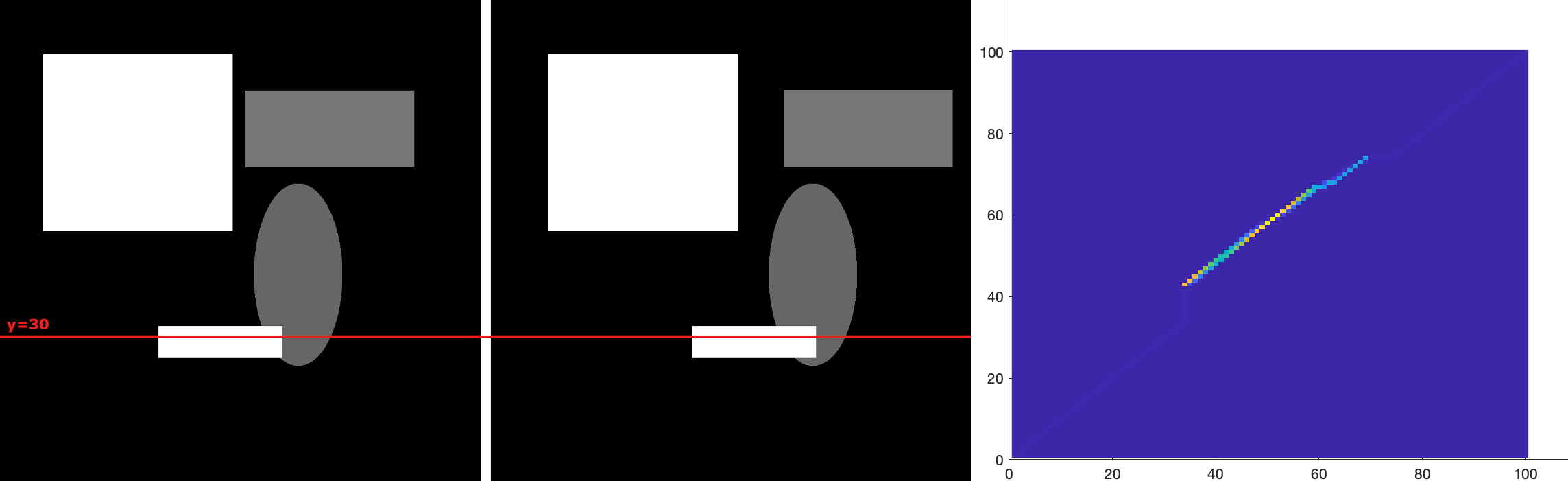}
\caption{In this figure we added a fourth object that has a superposition with the circle.
On the right we show the transport plan at $y=30$.}
\end{figure}

\begin{figure}[h]
\includegraphics[width=7cm]{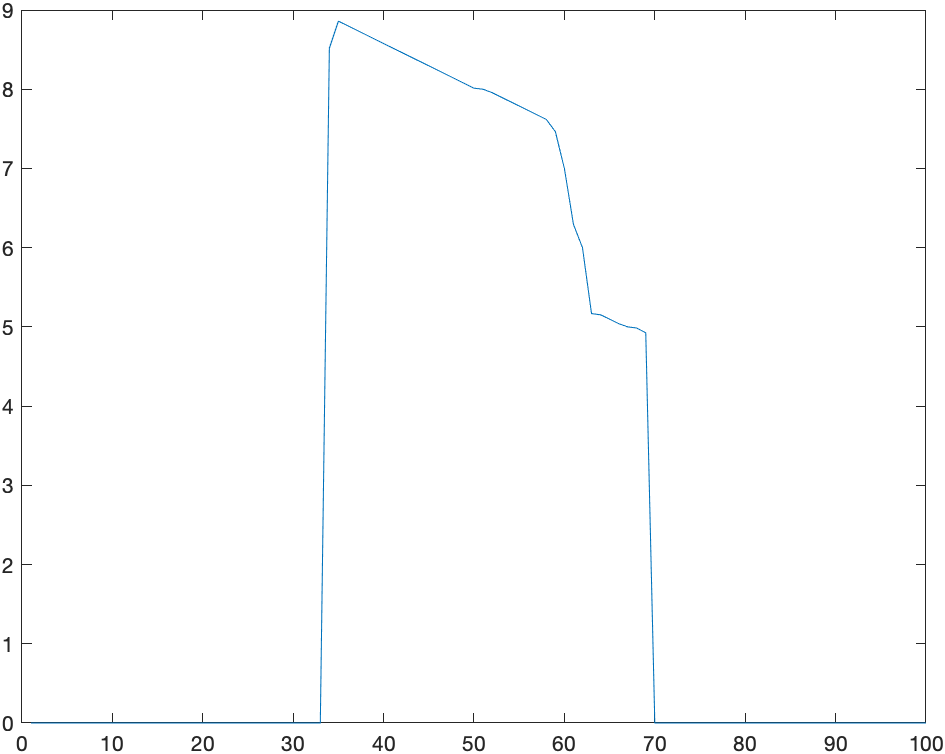}
\caption{The function obtained from the optimal plan at $y=30$. As we can see
the disparity shift for the leftmost object is $9$.}
\end{figure}

\begin{rmk}
We observe from the image below that the value of the compression coefficient
is not constant over the $X$ points of $[1,d]$. This is a consequence of the discretization
of the domain. If we want to evaluate $f_\gamma$ at $X$ we start by imposing, without
loss of generalities, that $i_0=1$. Then, we have two cases for any point $i$ in $X$.

If $\lfloor \varphi (i-1) \rfloor<\lfloor \varphi i \rfloor$, then
\begin{align*}
f_\gamma(i)&=\frac{1}{\varphi}\left(\{\varphi i\}\cdot \lceil \varphi i \rceil + (1-\{\varphi \cdot (i-1)\})\cdot \lfloor \varphi i\rfloor\right)-i\\
&=\left(1 -\frac{1}{\varphi}\right)\cdot \lfloor \varphi i \rfloor.
\end{align*}

While, if $\lfloor \varphi (i-1) \rfloor =\lfloor \varphi i\rfloor$, then
\begin{align*}
f_\gamma(i)&=\lceil \varphi i \rceil -i.
\end{align*}
The derivative operator in a discrete setting is defined as $\Delta f(i):=f(i+1)-f(i)$. We observe that
therefore $\Delta f_\gamma$ takes the same value for two successive entries $i, i+1$
only when this value is $1-\frac{1}{\varphi}$ and this allows to recover the $\varphi$ value
and to check our implementation.
\end{rmk}

\begin{figure}[h]
\includegraphics[width=9cm]{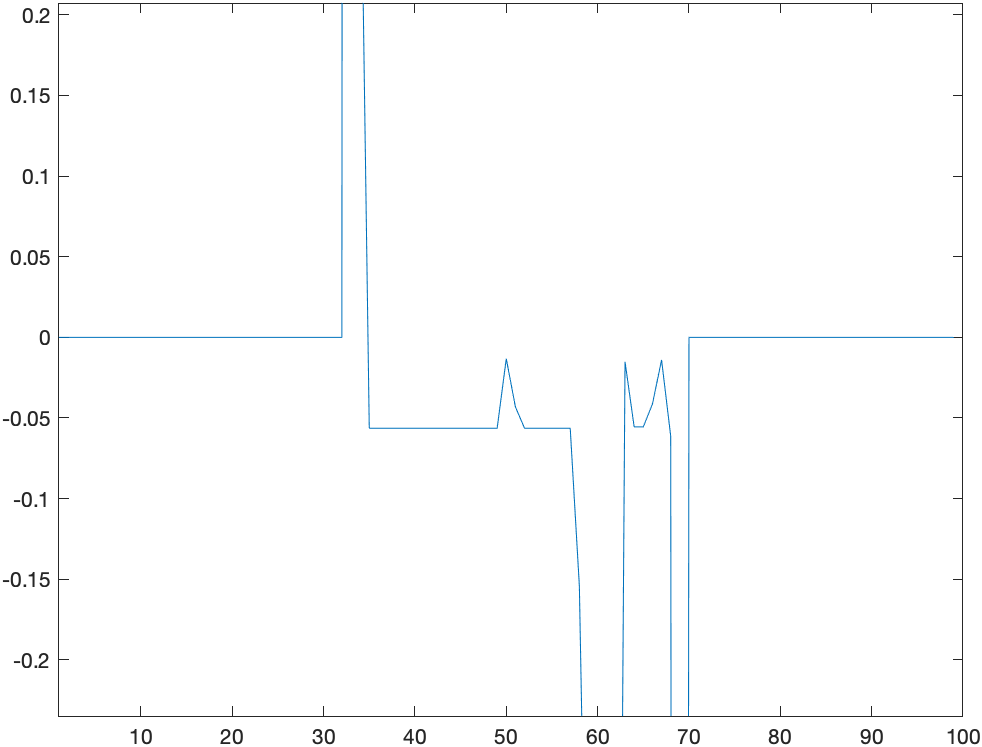}
\caption{The graph of the discrete derivative $\Delta f_\gamma$. We
recover from this that $1-\frac{1}{\varphi}\approx -0.05637$ and therefore
$\varphi\approx 0.95$, which coincides with the true estimate.}
\end{figure}




\bibliography{mybiblio}{}
\bibliographystyle{plain}

\end{document}